\def\zbar{\overline{z}}
\def\tr{{\rm tr}}

\def\IH{{\mathbb H}^3}

\def\IR{{\mathbb R}}

\def\ID{{\mathbb D}}

\def\IS{{\mathbb S}}
\def\IZ{{\mathbb Z}}
\def\IH{{\mathbb H}}
\def\IC{\mathbb C}

\def\oC{\hat{\IC}}

\def\zbar{{\overline{z}}}
\def\hbar{{\overline{h}}}
\def\zetabar{{\overline{\zeta}}}

\documentclass{article}

\usepackage{amsthm, amssymb}
\usepackage{amsmath}
\usepackage{amsfonts}
\usepackage{epsfig,multicol}
\usepackage{pstricks}
\usepackage{graphicx}

\title{{\bf Random Kleinian Groups}, {\bf I} \\
Random Fuchsian Groups.}

\author{Gaven Martin and Graeme O'Brien \thanks{Research supported in
part by grants from the N.Z.
Marsden Fund.  This work forms part of G. O'Brien's Thesis. \newline \newline AMS
(1991) Classification.
Primary 30C60, 30F40, 30D50, 20H10, 22E40, 53A35, 57N13, 57M60}  }

\date{}

\begin{document}

\maketitle
\newtheorem{theorem}{Theorem}[section]    
                                           
\newtheorem{lemma}[theorem]{Lemma}         
\newtheorem{corollary}[theorem]{Corollary} 
\newtheorem{remark}[theorem]{Remark}       
\newtheorem{definition}[theorem]{Definition}
\newtheorem{conjecture}[theorem]{Conjecture}
\newtheorem{proposition}[theorem]{Proposition}
\newtheorem{example}[theorem]{Example}

\newcommand{\param}{(\gamma,\beta,\beta')}
\newcommand{\parfour}{(\gamma,\beta,-4)}

\numberwithin{equation}{section}

\newcommand{\abs}[1]{\lvert#1\rvert}

\renewcommand{\theequation}{\thetheorem} 
                    
\makeatletter
\let \c@equation=\c@theorem

\begin{abstract}  We introduce a geometrically natural probability measure on the group of all M\"obius transformations of the circle.  Our aim is to study ``random'' groups of M\"obius transformations,  and in particular random two-generator groups.  By this we mean groups where the generators are selected randomly.  The probability measure in effect establishes an isomorphism between random $n$-generators groups and collections of $n$ random pairs of arcs on the circle.  Our aim is  to estimate the likely-hood that such a random group is discrete,   calculate the expectation of their associated parameters,  geometry and topology,  and to test the effectiveness of tests for discreteness such as J\o rgensen's inequality. 
\end{abstract}

\maketitle

\section{Introduction.}   In this paper we introduce the notion of a random Fuchsian group.  For us this will mean a finitely generated  Fuchsian group where the generators are selected from a geometrically natural probability measure on the space of M\"obius transformations of the circle.  Our ultimate aim is to study random Kleinian groups,  but the Fuchsian case is quite distinct in many ways - for instance the set of precompact cyclic subgroups (generated by elliptic elements) has nonempty interior in the Fuchsian case,  and therefore will have positive measure in any reasonable probability measure we might seek to impose.  Whereas for Kleinian groups this is not the case.  However,  the motivation for the probability measure we chose is similar in both cases.  We seek something ``geometrically natural'' and with which we can compute.  We should expect that almost surely (that is with probability one) a finitely generated subgroup of the M\"obius group is free.  We shall see that the probability a random two generator group is discrete is greater than $\frac{1}{20}$,  a value we conjecture to being close to optimal,  and this value is certainly less than $\frac{1}{4}$.  If we condition by choosing only hyperbolic elements,  this probability becomes $\frac{1}{5}$.  The cases where we condition by choosing two parabolic elements,  both in the Fuchsian and Kleinian cases is discussed in a sequel \cite{MOY} as rather more theory is required to get precise answers. Here we give a bound of $\frac{1}{6}$ in the Fuchsian case,  the actual value being approximately $0.3148\ldots$. 

\medskip

Here we also consider such things as the probability that the axes of hyperbolic generators cross.   This allows us to get some understanding of the likely-hood of different topologies arising. For instance if we choose two random hyperbolic elements with pairwise disjoint isometric circles,  the quotient space is either the two-sphere with three holes,  or a torus with one hole.  The latter occurring with probability $\frac{1}{3}$.   

\medskip

To study these questions of discreteness we set up a topological isomorphism between $n$ pairs of random arcs on the circle and $n$-generator Fuchsian groups.  We determine the statistics of a random cyclic group completely,  however, the statistics of commutators of generators is an important challenge with topological consequences and which we only partially resolve.

\section{Random Fuchsian Groups.}  We introduce specific definitions in the context of Fuchsian groups. These will naturally motivate more general definitions for the case of Kleinian groups in later work.

\medskip

If $A\in PSL(2,\IC)$ has the form
\begin{equation}\label{Fspace}
A =\pm \left( \begin{array}{cc} a & c \\ \bar c & \bar a \end{array} \right), \hskip15pt |a|^2-|c|^2 = 1,
\end{equation}
then the associated linear fractional transformation $f:\oC\to\oC$ defined by
\begin{equation}\label{fdef}
f(z) = \frac{az + c}{\bar c z + \bar a}
\end{equation}
preserves the unit circle since  $\left| \frac{az + c}{\bar c z + \bar a} \right| = |\zbar| \left| \frac{az + c}{\bar a \zbar+\bar c |z|^2} \right|$,  with the implication that $|z|=1$ implies $|f(z)|=1$.

The rotation subgroup ${\bf K}$ of the disk,  $z\mapsto \zeta^2 z$, $|\zeta|=1$,   and the nilpotent or parabolic subgroup ${\bf P}$ (conjugate to the translations) have the respective representations
\[  \left( \begin{array}{cc} \zeta & 0 \\ 0 & \bar \zeta \end{array} \right), \;\;\;|\zeta|=1, \hskip15pt   \left( \begin{array}{cc} 1+it& t \\ t & 1-i t\end{array} \right),\;\;\; t\in \IR . \]
The group of all matrices satisfying (\ref{Fspace}) will be denoted ${\cal F}$.  It is not difficult to construct an algebraic isomorphism ${\cal F}\equiv PSL(2,\IR)\equiv Isom^+(\IH^2)$,  the isometry group of two-dimensional hyperbolic space (see \cite{Beardon})   and we will often abuse notation by moving between $A$ and $f$ interchangeably.  Despite some efforts to use $PSL(2,\IR)$,  we feel the approach we take is geometrically more natural when working in ${\cal F}$. In particular,   our measures are obviously invariant under the action of the compact group ${\bf K}$.  We also seek distributions from which we can make explicit calculations and are geometrically natural (see in particular Lemma \ref{3.2}).

\medskip

We therefore impose the following distributions on the entries of this space of matrices ${\cal F}$.  We select 
\begin{itemize}
\item(i) $\zeta=a/|a|$ and $\eta=c/|c|$ are chosen uniformly in the circle $\IS$, with arclength measure,    and 
\item(ii) $t=|a|\geq 1$ is chosen so that
\[ 2\arcsin(1/t) \in [0,\pi] \]
is uniformly distributed.    
\end{itemize} 
Notice that the product  $\zeta\eta$ is uniformly distributed on the circle as a simple consequence of the rotational invariance of arclength measure.  Further,  this measure is equivalent to the uniform probability measure $\arg(a)\in [0,2\pi]$.  It is  thus clear that this selection process is invariant under the rotation subgroup of the circle.   
\medskip

Next,  if $\theta$ is uniformly distributed in $[0,\pi]$,  then the probability distribution function (henceforth p.d.f.) for $\sin \theta$ is $\frac{1}{\pi} \frac{1}{\sqrt{1-y^2}}$ for $y\in [-1,1]$.  Since $t\mapsto1/t$,  for $t>0$ is strictly decreasing,  we can use the change of variables formula for distribution functions to deduce the p.d.f. for $|a|$. 
\begin{lemma}  The random variable $|a|\in [1,\infty)$ has the p.d.f.
\[ F_{|a|}(x)= \frac{2}{\pi} \; \frac{1}{x\sqrt{x^2-1}} \]
\end{lemma}
Next notice that the equation $1+|c|^2=|a|^2$ tells us that $\arctan(\frac{1}{|c|})$ is also uniformly distributed in $[0,\pi]$.

\bigskip

Another equivalent formulation is the following.  We require that the matrix entries $a$ and $c$ have arguments $\arg(a)$ and $\arg(c)$ which are uniformly distributed on $\IR \mod 2\pi$.  We  write this as $\arg(a)\in_u [0,2\pi]_\IR$ and $\arg(c)\in_u [0,2\pi]_\IR$.  We illustrate with a lemma.

\begin{lemma}  If $\arg(a),\arg(b)\in_u [0,2\pi]_\IR$,  then $\arg(ab),  \arg(a/b) \in_u [0,2\pi]_\IR$.
\end{lemma}
\noindent {\bf Proof.}  The usual method of calculating probability distributions for combinations of random variables via characteristic functions shows that if $\theta,\eta$ are selected from a uniformly distributed probability measure on $[0,2\pi]$,  then the p.d.f. for $\theta+\eta\in [0,4\pi]$ is given by
\begin{equation}\label{f93i} g(\zeta)=\left\{ \begin{array}{lllll}
\frac{\zeta}{8\pi^2}&0\leq \zeta< 2\pi\\
\\
\frac{4\pi-\zeta}{8\pi^2}&2\pi \leq \zeta \leq 4\pi.
 \end{array}\right.\end{equation}
We reduce mod $2\pi$  and observe 
\[ \frac{\zeta}{8\pi^2} + \frac{4\pi-\zeta}{8\pi^2} =  \frac{1}{2\pi}  \] 
and this gives us once again the uniform probability density on $[0,2\pi]$. The result also follows for $a-b$ as clearly $-b \in_u [0,2\pi]_\IR$ and $a-b=a+(-b)$. \hfill $\Box$

\bigskip
 
 \begin{corollary}\label{cor2.6}  If $a\in_u [0,2\pi]_\IR$ and $k\in \IZ$,  then $ka \in_u [0,2\pi]_\IR$. 
\end{corollary}
 
 It now follows that any finite integral linear combination of variables $a_i\in_u [0,2\pi]_\IR$ has the same distribution.
 
 \medskip
 
 In what follows we will also need to consider variables supported in $[0,\pi]$ or smaller subintervals and as above we will write this as  $a\in_u [0,\pi]_\IR$ and so forth.
 
 \medskip
 
 In a moment we will calculate some distributions naturally associated with M\"obius transformations such as traces and translation lengths.  Every M\"obius transformation of the unit disk $\ID$ can be written in the form
 \begin{equation}\label{mob}
 z\mapsto \zeta^2 \,  \frac{z-w}{1-\bar w z}, \hskip10pt |\zeta|=1, w\in \ID
 \end{equation}
 Thus one could consider another approach by choosing distributions for $\zeta\in \IS$ and $w\in\ID$.   It seems clear one would want $\zeta$ uniformly distributed in $\IS$.  The real question is by what probability measure should $w$ be chosen on $\ID$ ?  If $w$ is chosen rotationally invariant,  then the choice boils down to probability measures on radii.  The choices we have made turn out as follows.  The matrix representation of (\ref{mob}) in the form (\ref{Fspace}) is 
 \[ \zeta^2 \,  \frac{z-w}{1-\bar w z} \leftrightarrow \left(\begin{array}{cc} \frac{\zeta}{\sqrt{1-|w|^2} } &  -\frac{\zeta w}{\sqrt{1-|w|^2}} \\- \frac{\zetabar \bar w}{\sqrt{1-|w|^2} }& \frac{\zetabar}{\sqrt{1-|w|^2}  }
   \end{array} \right) \]
Hence $\zeta$ and $w/|w|$ will be uniformly distributed in $\IS$.  Then,  $|w|<1$ necessarily and \[  \arccos(|w|)=\arcsin(\sqrt{1-|w|^2} ) \in [0,\pi/2] \]
is uniformly distributed and we find $|w| = |f(0)|$ has the p.d.f. $\frac{2}{\pi  \sqrt{1-y^2}}$,  $y\in [0,1])$.

\begin{corollary}  Let $f\in {\cal F}$ be a random M\"obius transformation.  Then the p.d.f. for $y=|f(0)|$ is $\frac{2}{\pi \sqrt{1-y^2}}$.  The expected value of $|f(0)|$ is 
\[ E[\;|f(0)|\; ]=\frac{2}{\pi} \int_{0}^{1} \frac{ y}{\sqrt{1-y^2}}\, dy = \frac{2}{\pi} = 0.63662\ldots \]
\end{corollary}
\noindent The hyperbolic distance here between $0$ and $f(0)$ is $\log \frac{1+|f(0)|}{1-|f(0)|} = \log \frac{\pi+2}{\pi-2} = 1.50494 \ldots $.
  
\section{Fixed points}  The fixed points of a random $f\in {\cal F}$ are solutions to the same quadratic equation and one should therefore expect some correlation.
From (\ref{fdef}) we see the fixed points are  the solutions to $az+c=z(\bar c z+ \bar a)$.  That is
\begin{equation}\label{fp}
z_\pm =\frac{1}{\bar c}\left( i \Im m(a) \pm \sqrt{\Re e(a)^2-1}\right), \hskip15pt |a|^2=1+|c|^2.
\end{equation}
We consider two cases and will soon establish that ${\rm Pr}\{|\Re e(a)|\leq 1\}=\frac{1}{2}$ so each case occurs with equal probability. 

\medskip

{\bf Case 1. } $f$ elliptic or parabolic.  Then $|\Re e(a)|\leq 1$ and so $\arg(z_\pm) = \frac{\pi}{2}+\arg(c)$.  Thus  the argument of both fixed points is the same and that angle is uniformly distributed in $[0,\pi]$.

\medskip

{\bf Case 2. } $f$ hyperbolic.  Then $\Re e(a) >1$ and $|z_\pm|=1$.  We calculate that the derivative 
\[ |f'(z_\pm)| = \frac{1}{|\bar c z_\pm +\bar a|^2}  =  \frac{1}{|i \Im m(a) \pm \sqrt{\Re e(a)^2-1}+\bar a|^2} = \frac{1}{|\Re e(a) \pm \sqrt{\Re e(a)^2-1}|^2}  \]
Hence $|f'(z_+)|<1$ and $z_+$ is an {\em attracting} fixed point,  with $z_-$ being {\em repelling}.

We have chosen $\arg(c)$ to be uniformly distributed and so the argument of either fixed point, say $z_+$, is uniformly distributed.  The interesting question is the distribution of the angle (at $0$) between the fixed points.   That is the argument of $z_+\overline{z_-}$.  This will reflect the correlation we are looking for.  This angle is easily seen to be the angle $\phi\in [0,\pi]$ where $\cos(\phi/2)= \Im m(a)/|c|$.  Then
\begin{eqnarray*}
\cos(\phi/2) & = & \Im m(a)/|c| = \frac{|a| \sin \theta}{\sqrt{|a|^2-1}} =   \frac{\sin \theta}{\cos \alpha}
\end{eqnarray*}
where we are able to assume that both $\theta$ and $\alpha$ are uniformly distributed in $[0,\pi/2]$ and we are conditioned by
$ {\sin \theta}\leq {\cos \alpha}.$

We will calculate the distribution of ${\sin \theta}/{\cos \alpha}$ carefully when we come to the calculation of the parameters determining a M\"obius group.  We report the p.d.f here as follows.

\begin{theorem}\label{p.d.f.1} The distribution of the random variable
\[ X=\frac{\sin(\theta)}{\cos(\alpha)} ,  \]
 for $\theta$ and $\alpha$ uniformly distributed in $[0,\pi/2]$ is given by the formula
\begin{equation}\label{dist1}
h_X(x)=\frac{4}{\pi^2x}\; 
\log{\frac{1+x}{1-x}},  \hskip10pt 0\leq x< 1.
 \end{equation}
\end{theorem} 
We can now use the change of variables formula to compute the p.d.f for $\phi/2$.  That is we want the distribution for $Y=\cos^{-1}(h_X(x))$,  given $h_X(x)\leq 1$.  We can compute this distribution to be
\[ h_Y(y) =  \frac{4}{\pi^2}\; \tan(y)\; 
\log{\frac{1+\cos(y)}{1-\cos(y)}}  \]
 \begin{theorem}\label{p.d.f.2} Let $\phi \in [0,\pi]$ be the angle subtended at $0$ by the fixed points of a random hyperbolic element in ${\cal F}$. 
Then the p.d.f. for $\eta = \phi/2$ is given by
\begin{equation}\label{dist2}
H_Y(\eta) =  \frac{4}{\pi^2}\; \tan(\eta)\; 
\log{\frac{1+\cos(\eta)}{1-\cos(\eta)}}   
 \end{equation}
\end{theorem} 
Some hyperbolic trigonometry reveals the the hyperbolic line between a pair of points $z_\pm\in \IS$ meets the closed disk of hyperbolic radius $r$ (denoted $\ID_\rho(r)$) when the angle $\phi$ formed at $0$ satisfies
\begin{equation}
\cosh(r) \geq \frac{1}{\sin(\phi/2)}.
\end{equation}
If $z_\pm$ are the fixed points of a hyperbolic element $f$,  then this hyperbolic line joining them is called the axis of $f$,  denoted ${\rm axis}(f)$.  We can therefore compute the probability that the axis of a random hyperbolic element meets $\ID_\rho(r)$ by setting $\delta=\sin^{-1}(1/\cosh(r) )$ and computing
\begin{eqnarray*}
{\cal P}({\rm axis}(f)\cap \ID_{\rho}(r) \neq \emptyset ) & = &  \frac{4}{\pi^2} \; \int_{0}^{\delta} \tan(\eta)\;  \log{\frac{1+\cos(\eta)}{1-\cos(\eta)}}   \; d\eta \\
& = &   \frac{4}{\pi^2} \; \int_{0}^{\tanh(r)} \frac{1}{x} \;  \log{\frac{1+x}{1-x}}   \; dx \\
& = & \frac{4}{\pi^2}  \big[ \text{Li}_2(\tanh (r))-\text{Li}_2(-\tanh (r))\big]
\end{eqnarray*}
Here $\text{Li}_2(s) = \sum_{1}^{\infty} n^{-2} s^n $ is a polylog function.   Thus,  for instance,  this probability exceeds $\frac{1}{2}$ as soon as $r>0.678\ldots$ and exceeds $0.95$ as soon as $r>2.24419$.

 \begin{center}
\scalebox{0.6}{\includegraphics*[viewport=-10 370 660 730]{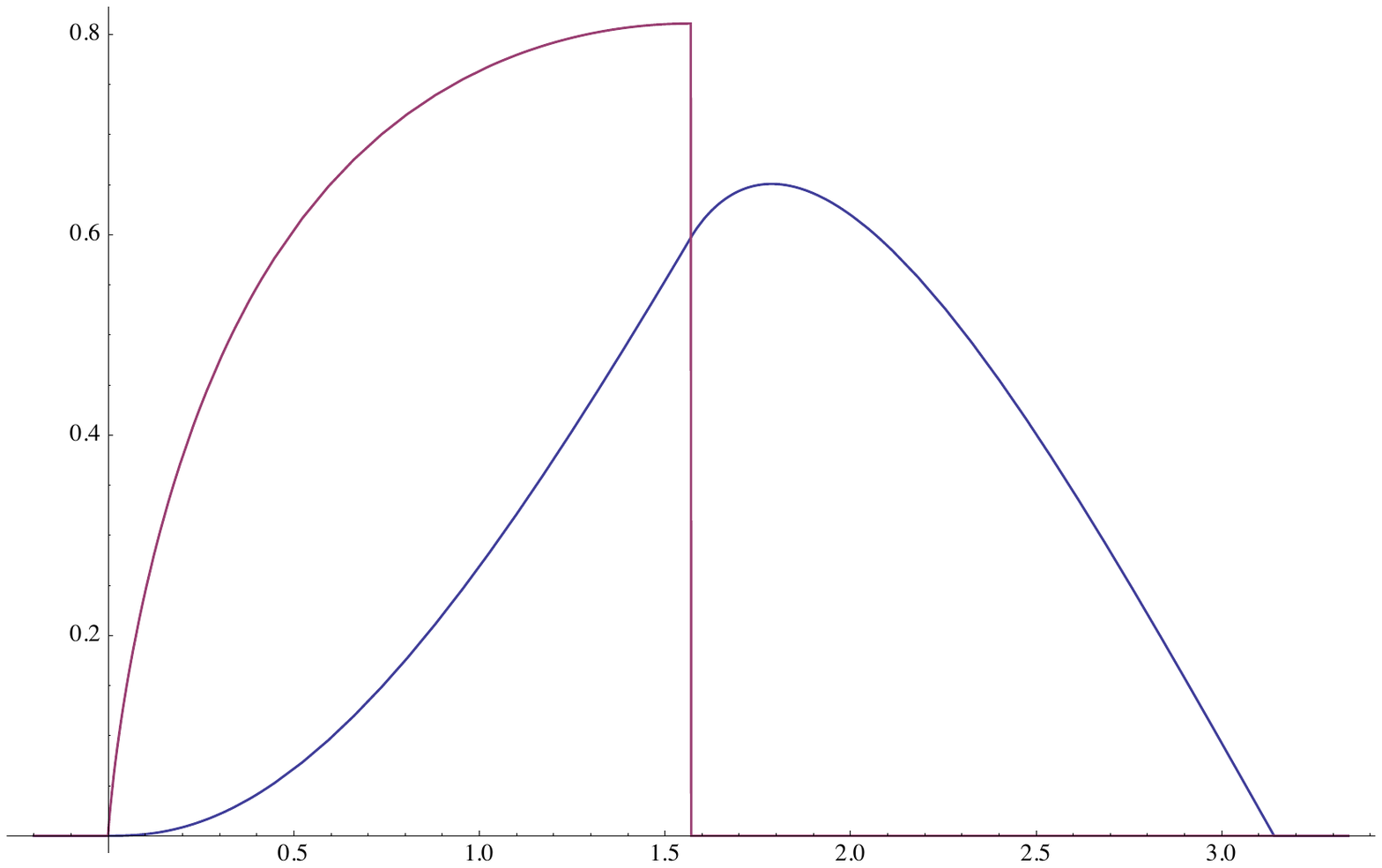}} \\
The p.d.f. $H_Y$ for the angle  $\phi/2$ between fixed points of a random hyperbolic  $f\in {\cal F}$ and the convolution $H_Y*H_Y$.
\end{center}

Now,  the bisector $\zeta_f$ of the smaller circular arc between the fixed points of a random hyperbolic element of $f$ is uniformly distributed on the circle.  Then,  given $f$ and $g$ random hyperbolic elements of ${\cal F}$ and angles $\phi_f$ and $\phi_g$ between their fixed points.  The p.d.f. for $\phi_f/2+\phi_g/2$ is the convolution $H_Y*H_Y$. We note that $e^{i\theta}= \xi = \zeta_f\overline{\zeta_g}$ is uniformly distributed as well.  Given $\xi$ the fixed points of $f$ and of $g$ intertwine (so that the axes cross) if both $\phi_f+\phi_g \geq 2\theta$ and $|\phi_f-\phi_g|<2\theta$.  We can use the distributions above to calculate these probabilities, but it is quite complicated and we will find another route to this probability a bit later.
 
\section{Isometric Circles and Traces.}
The isometric circles of the M\"obius transformation $f$ defined at (\ref{fdef}) are defined to be the two circles
\[  C_+ =\Big\{|z+\frac{\bar a}{\bar c} | = \frac{1}{|c|}\Big\},  \hskip10pt C_-=\Big\{z:|z-\frac{a}{\bar c}|=\frac{1}{|c|} \Big\} \]
which are paired by the action of $f$ and $f^{-1}$,   $f^{\pm1}(C_{\pm})=C_{\mp}$.  The {\em isometric disks} are the finite regions bounded by these two circles.  

Since $|a|^2=1+|c|^2\geq 1 $,  both these circles meet the unit circle in an arc of angle $\theta\in [0,\pi]$.  Some elementary trigonometry reveals that
\begin{equation}\label{a}
\sin \frac{\theta}{2} = \frac{1}{|a|}
\end{equation}
Thus by our choice of distribution for $|a|$ we obtain the following key result.

\begin{lemma} \label{3.2} The arcs determined by the intersections of the finite disks bounded by the the isometric circles of $f$,  where $f$ is chosen according to the distribution (i) and (ii),  are centred on uniformly distributed points of $\IS$ and have arc length uniformly distributed in $[0,\pi]$.  
\end{lemma}
It is this lemma which supports our claim that the p.d.f. on ${\cal F}$ is natural and suggests the way forward for an analysis of random Kleinian groups.

\medskip 

The isometric circles of $f$ are disjoint if 
$\left|\frac{ a}{\bar c}+\frac{\bar a}{\bar c} \right| \geq \frac{2}{|c|} $.
This occurs if 
\[ |\tr(f)| = |a+\bar a| = 2 |\Re e(a)| \geq 2\]
Since the  disjointness of isometric circles has important geometric consequences  we will need  to find the p.d.f. for  the random variable $t=|\tr(f)|$. 

\medskip

As $|\Re e(a)| = |a||\cos(\theta)|$,  for a fixed $\theta\in [0,\pi/2]$,  the probability
\begin{eqnarray}\label{tracedistn}
Pr[ \{ |a|\geq 1/\cos \theta \}] =1 - \frac{2}{\pi} \int_{1}^{1/\cos \theta} \frac{dx}{x\sqrt{x^2-1}} = 1-\frac{2}{\pi}\theta
\end{eqnarray}
As $a/|a|$ is unformly distributed on the circle,  we have $\theta|[0,\pi/2]$ uniformly distributed in $[0,\pi/2]$.  Therefore using the obvious symmetries we may calculate that
\begin{eqnarray*}
Pr[ \{ |a+\bar a| \geq 2 \}] = \frac{2}{\pi} \; \int_{0}^{\pi/2} \; 1-\frac{2}{\pi}\theta \; d\theta = \frac{1}{2}.
\end{eqnarray*}
\begin{corollary}  Let $f\in {\cal F}$ be a M\"obius transformation chosen randomly from the distribution described in (i) and (ii).  Then the probability that the isometric circles of $f$ are disjoint is equal to $\frac{1}{2}$.  
\end{corollary}
Therefore we have the following simple consequence concerning random cyclic groups.
\begin{corollary}  Let $f\in {\cal F}$ be a M\"obius transformation chosen randomly from the distribution described in (i) and (ii).  Then the probability that the cyclic group $\langle f \rangle$ is discrete is equal to $\frac{1}{2}$.  
\end{corollary}
\noindent{\bf Proof.}  The matrix $A\in SL(2,\IC)$ represents the hyperbolic M\"obius transformation $f$ if and only if $-2 \leq \tr A \leq 2$.  This occurs with probability $\frac{1}{2}$.  The matrix $A$ represents an elliptic transformation of finite order,  or a parabolic transformation if and only if $\tr(A) = \pm 2\cos(p\pi/q)$,  $p,q\in \IZ$,  and this set is countable and therefore has measure zero.  The result follows. \hfill $\Box$.
\medskip

We now note the following trivial consequence.
\begin{corollary}  Let $f,g\in {\cal F}$ be M\"obius transformations chosen randomly from the distribution described in (i) and (ii).  Then the probability that the group $\langle f,g \rangle$ is discrete is no more than $\frac{1}{4}$.  
\end{corollary}

\bigskip

Actually we can use (\ref{tracedistn}) to determine the p.d.f. for $|\tr(A)|$.  We will do this two ways.  

First, for $s\geq 2$, 
\begin{eqnarray*}\label{tracedist}
Pr[ \{ |\tr(A)| \geq s \}]  & = & Pr[ \{ 2|a|\cos \theta \geq s \}] = Pr[ \{ |a|\geq s/(2\cos \theta) \}] \\
& = & 1 - \frac{4}{\pi^2} \int_{0}^{\pi/2}    \; \int_{1}^{s/2\cos \theta} \frac{dx}{x\sqrt{x^2-1}} \; d\theta \\
& = & 1 - \frac{4}{\pi^2} \int_{0}^{\pi/2}    \cos^{-1} \Big( \frac{2\cos \theta}{s} \Big)\; d\theta 
\end{eqnarray*}
We can now differentiate this function of $s$ under the integral,  integrate with respect to $\theta$ (using the symmetry to reduce it to being over $[0,\pi/2]$), to obtain the probability density function for $|\tr(A)|$ (for $|\tr(A)|\geq 2$),
\begin{eqnarray} 
 F[s] = \frac{4}{\pi^2\, s} \cosh^{-1}\Big( \frac{s}{\sqrt{s^2-4}} \Big)  , \hskip10pt s\geq 2.
 \end{eqnarray}
 This gives the distribution for $\tr^2 A$ as 
\begin{eqnarray} 
G[t] = \frac{2}{\pi^2\, t } \cosh^{-1}\Big( \frac{\sqrt{t}}{\sqrt{t-4}} \Big)  = \frac{2}{\pi^2 \, t} \log \frac{\sqrt{t}+2}{\sqrt{t-4}} , \hskip10pt t\geq 4.
 \end{eqnarray}
 Then the random variable $\beta = \tr^2A-4\geq 0$ has distribution
\begin{eqnarray} \label{beta1}
G[\beta] =  \frac{1}{\pi^2 (\beta+4)} \log \left( 1+\frac{8+4\sqrt{\beta+4}}{\beta} \right), \hskip10pt \beta \geq 0.
 \end{eqnarray}
 We could now follow through a similar,  but more difficult, calculation to determine the distribution for $\beta$ in the interval $-4 \leq  \beta \leq 0$.  It turns out to be 
 \begin{eqnarray} 
G[\beta] = \frac{1}{\pi^2(\beta+4) }\log \left( \frac{2+\sqrt{\beta + 4}}{2 - \sqrt{\beta + 4}}  \right), \hskip15pt \beta\in [-4,0].
\end{eqnarray} 
We will return to this in a moment through a different approach as we can immediately use (\ref{beta1})
to find the distribution of the translation length of hyperbolic elements.   

As we have seen,  every element $f\in{\cal F}$ which is not elliptic (conjugate to a rotation,  equivalently $\beta(f)\in [-4,0)$) or parabolic (conjugate to a translation,  equivalently $\beta(f)=0$) fixes two points on the circle and the hyperbolic line ${\rm axis}(f)$ with those points as endpoints.  The transformation acts as a translation by constant hyperbolic distance $\tau(f)$ along its axis.  This number $\tau(f)$  is called the {\em translation length} and is related to the trace via the formula \cite{GehMar}
\[ \beta(f) = 4 \sinh^2 \;  \frac{\tau}{2},  \hskip10pt \tau =  \cosh^{-1}\left(1+\frac{\beta}{2}\right)\]
We obtain the distribution for $\tau = \tau(f)$ from the change of variables formula for p.d.f.  using (\ref{beta1})
\begin{eqnarray*}
H[\tau]  & = &   \frac{2}{\pi^2} \,  \tanh  \frac{\tau}{2}  \;\log \left(\frac{ \cosh \,  \frac{\tau}{2} +1}{ \cosh\, \frac{\tau}{2} -1} \right)  \\
& = &  - \frac{4}{\pi^2} \; \tanh  \frac{\tau}{2} \; \log \tanh\,  \frac{\tau}{4}  \\
\end{eqnarray*}
Unlike our earlier distribution $G$,  the p.d.f for $\tau$ has all moments.  In particular once we observe
\[ \int_{0}^{\infty} t \,\tanh \,\frac{t}{2} \log\Big[\tanh \, \frac{t}{4}\Big] dt = -\pi^2 \log 2 \]
we have the following theorem.
\begin{theorem}  For randomly selected hyperbolic $f\in {\cal F}$ the p.d.f. for the translation length $\tau=\tau(f)$ is 
\begin{equation}
H[\tau]  = - \frac{4}{\pi^2} \; \tanh  \frac{\tau}{2} \; \log \tanh\,  \frac{\tau}{4}
\end{equation}
(illustrated below) and the expected value of the translation length is 
\begin{equation}
E[[\tau]] = 4\log 2 \approx 2.77259 \ldots
\end{equation}
\end{theorem}
 \begin{center}
\scalebox{0.5}{\includegraphics*[viewport=-100 470 560 780]{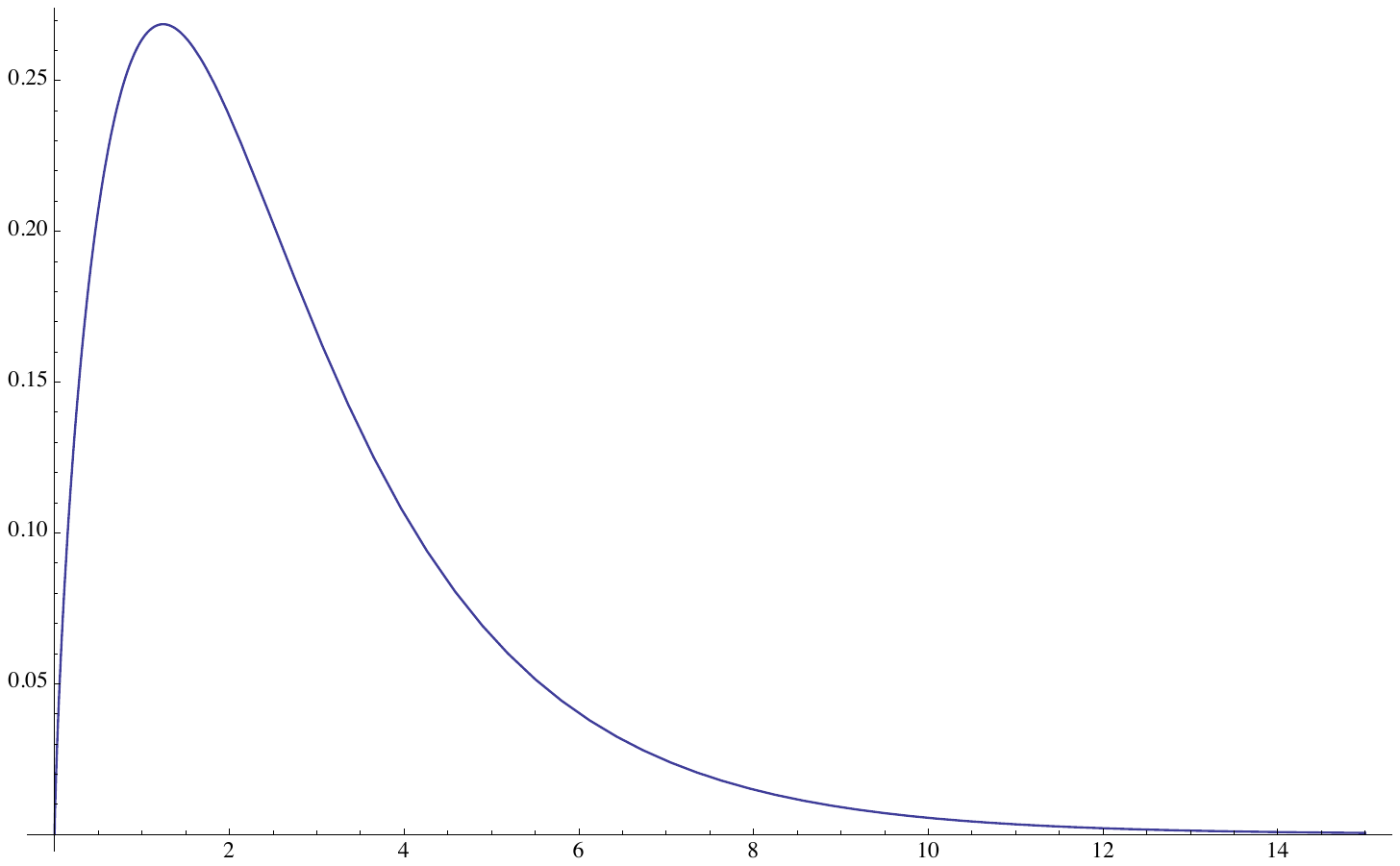}} \\
The p.d.f for the translation length $\tau$ of a random hyperbolic element of ${\cal F}$.
\end{center}

However there is another way to see  these results and which is more useful in what is to follow in that it more clearly relates to the geometry. 

\section {The parameter $\beta=\tr^2(A)-4$}

We being with the following theorem.

\begin{theorem}\label{th_beta} If a M\"obius transformation $f$ is randomly chosen in ${\cal F}$,   then 
\begin{equation} \beta(f)=4\left(\frac{\cos^2(\theta)}{\sin^2(\alpha)}-1\right)\;\;\;\;\;\;\;\;\theta\in_u[0,2\pi],\;\alpha\in_u\left[0,\frac{\pi}{2}\right]\label{beta1*}
\end{equation}
where $2\alpha$ is the arc length intersection of the isometric circles of $f$ with $\IS$ and $\theta$ is the argument of the leading entry of  $A$,  the matrix representative for $f$.
\end{theorem}
\noindent{\bf Proof.}  Let  $A=\left(\begin{array}{cc}  a & c \\ \bar c & \bar a \end{array}\right)$.  Then 
\[ \beta = \tr^2 A -4 = [2 \Re e(a)]^2-4 = 4 |a|^2 \cos^2(\theta)-4 \]
and the result follows by (\ref{a}) and Lemma \ref{3.2}. \hfill $\Box$

\bigskip

\begin{theorem}\label{th_c2_s2} The distribution of the random variable
\[ w=\frac{\cos^2(\theta)}{\sin^2(\alpha)} ,  \hskip10pt{\rm for} \;\;\;\; \theta\in_u[0,2\pi] \;\;\; {\rm and} \;\;\; \alpha\in_u[0,\frac{\pi}{2}]\]
 is given by the formula

\begin{equation}\label{dist1}
h(w)=\frac{1}{\pi^2w} 
\log \left| \frac{\sqrt{w}+1}{\sqrt{w}-1}\right|,  \hskip20pt w\geq 0.
\end{equation}\end{theorem}
\noindent {\bf Proof.}
The p.d.f's of $x=\cos^2(\theta)$ and $y=\sin^2(\alpha)$ are 
\begin{equation}\label{c2f2} f(x) =\frac{1}{\pi\sqrt{x(1-x)}}\;{\rm for} \;\cos^2(\theta),  \;\; {\rm and} \;\;
f(y)=\frac{1}{\pi\sqrt{y(1-y)}}\; {\rm for} \;\sin^2(\alpha).
\end{equation}
 and these are identically distributed when both $\theta$ and $\alpha$ are identically distributed.  They are also monotonic for $x,y\in[0,\frac{1}{2})$ and also for $x,y\in(\frac{1}{2},1]$ and as the distributions are anti-symmetric about $\frac{1}{2}$.  Therefore we can use the change of variables formula and the Mellin convolution to compute the p.d.f.  Write $x=\cos^2(\theta)$, $y=\sin^2(\alpha)$ and $w=\frac{\cos^2(\theta)}{\sin^2(\alpha)}$.  We use the Mellin convolution for quotients as in \cite{Springer}, noting that the distributions $f(x)$ and $f(y)$ in are identical.  For $x,y\in(0,1)$ the upper integration limits for the convolution integrals  will be $y<1\times\frac{1}{w}$ whenever $w> 1$ and $y<1$ otherwise, accordingly the Mellin convolution for the quotient of the p.d.f's over $(0,\infty)$ is calculated as follows where  we have ensured the piecewise differentiability of the integrand.
\begin{equation}\label{hw1}
h(w)=\left\{\begin{array}{ll}
\int_{0}^{1} y\;f(x)f(y) dy\;\;\;\;\;\;&w<1\\
\\
\int_0^{\frac{1}{w}} y\;f(x)f(y) dy&w>1
\end{array}\right.
\end{equation}
and the indefinite integral embedded in both components of \eqref{hw1} is given as
\begin{eqnarray} \label{p.d.f_c}
\int y\;f(yw)f(y)dy&= & \int y\;\frac{1}{\pi\sqrt{yw(1-yw))}}\frac{1}{\pi\sqrt{y(1-y)}}dy  \nonumber \\
& = & \frac{1}{\pi^2 \sqrt{w}}\int\frac{1}{\sqrt{(1-y)(1-yw)}}dy \nonumber \\
&=& \frac{2}{\pi^2w} \; \log\left(w\sqrt{(y-1)}+\sqrt{w(yw-1})\right) .
\end{eqnarray} 
Simplification of the $\log$ term in (\ref{p.d.f_c}) yields 
\begin{eqnarray*}
\lefteqn{\log\left(w\left(w(2y-1)-1+2\sqrt{w(y-1)(yw-1)}\right)\right)} \\
& = & 
 \left\{ \begin{array}{lll}
e_0&=\log(-w(w+1-2\sqrt{w}))\;\;\;\;&{\rm at} \;y=0\\
e_1&=\log(w(w-1))&{\rm at} \;y=1\\
e_{\frac{1}{w}}&=\log(-w(w-1))&{\rm at} \;y=\frac{1}{w}.
\end{array}\right.\end{eqnarray*} 
and accordingly the definite integrals in (\ref{hw1}) evaluate to
\[  
\int_0^1 y\;f(yw)f(y)dy = \frac{1}{\pi^2w}(e_1-e_0), \;\;\; 
\int_0^{\frac{1}{w}} y\;f(yw)f(y)dy= \frac{1}{\pi^2w}(e_{1/w}-e_0).
\]
If we now let $v=\sqrt{w}$,  then 
\begin{eqnarray*}
e_1-e_0&=& \log(w(w-1))-\log(-w(w+1-2\sqrt{w}))=\log\left(\frac{w(w-1)}{-w(w+1-2\sqrt{w}}\right))\\
&=&\log\left(\frac{(v^2-1)}{-(v^2+1-2v}\right)=\log\left(\frac{(v-1)(v+1)}{-(v-1)^2}\right)= \log\left(\frac{1+\sqrt{w}}{1-\sqrt{w}}\right)
\end{eqnarray*}
and
\begin{eqnarray*}
e_{1/w}-e_0&=& \log(-w(w-1))-\log(-w(w+1-2\sqrt{w}))\\ & = & \log\left[\frac{\log(-w(w-1))}{-w(w+1-2\sqrt{w}}\right]\\
&=& \log\left(\frac{-(v^2-1)}{-(v^2+1-2v}\right)= \log\left(\frac{\sqrt{w}+1}{\sqrt{w}-1}\right)
\end{eqnarray*} 
We therefore deduce that the distribution of $w=\frac{\cos^2(\theta) }{\sin^2(\alpha)}$ is given by (\ref{dist1}) as claimed. \hfill $\Box$
 
 \bigskip
 
 From this,  and a little obvious manipulation to see these formulas actually agree with those obtained earlier,  we obtain the result we were looking for.

\begin{theorem}\label{p.d.f._beta} The distribution of $\beta(f)$ for $f$ randomly chosen from ${\cal F}$ is given by
\begin{equation}G[\beta] =\frac{4}{\pi^2(\beta+4)} \; 
\log \Big|\frac{\sqrt{\beta+4}+2}{\sqrt{\beta+4}-2} \Big|, \hskip10pt \beta\geq -4 \end{equation}\end{theorem}
 \begin{center}
\scalebox{0.4}{\includegraphics{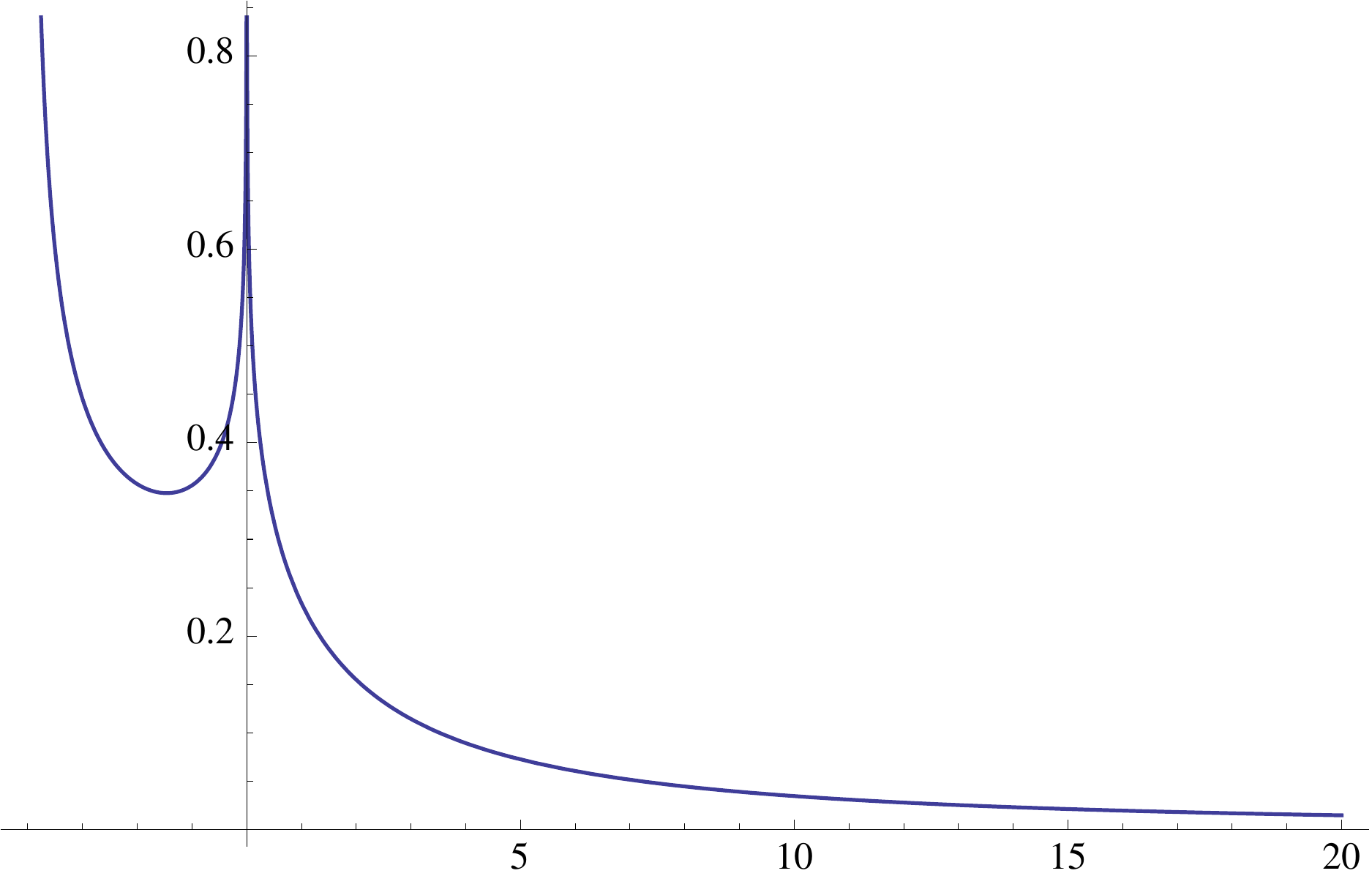}} \\
The p.d.f for the parameter $\beta(f)$ for a random element $f\in {\cal F}$.
\end{center}

\section{The topology of the quotient space.}

Topologically there are two surfaces whose fundamental group is isomorphic to $F_2$,  the free group on two generators.  These are the $2$-sphere with three holes $\IS^{2}_{3}$,  and the Torus with one hole $T^{2}_{1}$.  Thus we can expect that a group $\Gamma=\langle f,g \rangle$ generated by two random hyperbolic elements of ${\cal F}$ if discrete,  has quotient  space
\[ \ID^2/\Gamma \in\{ \IS^{2}_{3},  T^{2}_{1} \} \]
We would like to understand the likely-hood of one of these topologies over the other.   The topology is determined by whether the axes of $f$ and $g$ cross (giving $ T^{2}_{1} $) or not (giving $\IS^{2}_{3}$).  This is the same thing as asking if the hyperbolic lines between the fixed points of $f$ and the fixed points of $g$ cross or not,  and this in turn is determined by a suitable cross ratio of the fixed points.  In fact,  the geometry of the commutator $\gamma(f,g)=\tr [f,g]-2$ determines not only the topology of the quotient,  but also the hyperbolic length of the shortest geodesic - it is represented by either $f$, $g$ or $[f,g]=fgf^{-1}g^{-1}$ and their Nielsen equivalents. In fact the three numbers $\beta(f),\beta(g)$ and $\gamma(f,g)$ determine the group $\langle f,g \rangle$ uniquely up to conjugacy.  Since we have already determined the natural probability densities for $\beta(f)$ and $\beta(g)$ we  need only identify the p.d.f. for $\gamma=\gamma(f,g)$ to find a conjugacy invariant way to identify random discrete groups.  Unfortunately this is not so straightforward and we do not know this distribution.  However important aspects of this distribution can be determined.

\subsection{Commutators and cross ratios.}

We follow Beardon \cite{Beardon} and define the cross ratio of four points $z_1,z_2,z_3,z_4\in \IC$ to be
\begin{equation}
[z_1,z_2,z_3,z_4] = \frac{(z_1-z_3)(z_2-z_4)}{(z_1-z_2)(z_3-z_4)}
\end{equation}
In order to address the distribution of $\gamma(f,g) = \tr [f,g]-2$  we need to understand the cross ratio distribution.  This is because of the following result from  \S 7.23 \& \S 7.24 \cite{Beardon} together with a little manipulation.
\begin{theorem}  Let $\ell_1$,  with endpoints $z_1,z_2$,  and $\ell_2$,  with endpoints $w_1,w_2$,  be hyperbolic lines in the unit disk model of hyperbolic space.  So $z_1,z_2,w_1,w_2\in \IS$,  the circle at infinity.  Let $\delta$ be the hyperbolic distance between $\ell_1$ and $\ell_2$,  and should they cross,  let $\theta\in [0,\pi/2]$ be the angle at the intersection.  Then
\begin{equation} \label{cxdist}
\sinh^2 \Big[ \frac{1}{2}(\delta+i\theta)\Big] \times   [z_1,w_1,z_2,w_2]   = -1
\end{equation}
\end{theorem}
The number $\delta+i\theta$ is called the {\em complex distance} between the lines $\ell_1$ and $\ell_2$ where we put $\theta=0$ if the lines do not meet.  The proof of this theorem is simply to use M\"obius invariance of the cross ratio and the two different models of the hyperbolic plane.   If the two lines do not intersect,  we choose the M\"obius transformation which sends the disk to the upper half-plane and $\{z_1,z_2\}$ to $\{-1,+1\}$ and $\{w_1,w_2\}$ to $\{-s,s\}$ for some $s>1$.   Then $\delta=\log s$ and 
\[ [-1,-s,1,s]=\frac{-4s}{(1-s)^2 } =  \frac{-4}{(e^{\delta/2}-e^{-\delta/2})^2 } = - \frac{1}{\sinh^2(\delta/2)} \]
 while if the axes meet at a finite point,  we choose a M\"obius transformation of the disk so the line endpoints are $\pm 1$ and $e^{\pm i\theta}$ and the result follows similarly.
 
 \medskip
 
 We next recall Lemma 4.2 of \cite{GehMar} which relates the parameters and cross ratios.
 
 \begin{theorem}  Let $f$ and $g$ be M\"obius transformations and let $\delta+i\theta$ be the complex distance between their axes.  Then
 \begin{equation}
 4 \gamma (f,g) = \beta(f) \, \beta(g) \, \sinh^2(\delta+i\theta).
 \end{equation}
 \end{theorem}
 We note from (\ref{cxdist}) that 
 \[ \sinh^2(\delta+i\theta) = \left(1-\frac{2}{[z_1,w_1,z_2,w_2]}\right)^2-1 \]
 For a pair of hyperbolics $f$ and $g$ we have $\beta(f),\beta(g)\geq 0$ with $\delta=0$ if the axes meet.  Thus the axes cross if and only if $\gamma<0$,  or equivalently 
\begin{equation} \label{cr1} [z_1,w_1,z_2,w_2] >  1.
\end{equation}
Actually to see the latter point,  we choose the M\"obius transformation which sends $z_1\mapsto 0$, $z_2\mapsto\infty$, $w_1\mapsto 1$.  Then $z_2\mapsto z$,  say,  and
\[  [z_1,w_1,z_2,w_2] = \frac{(0-1)(\infty-z)}{(0-\infty)(1-z)}=\frac{1}{1-z} \] 
The image of the axes (and therefore the axes themselves) cross when $z<0$,  equivalently when  (\ref{cr1}) holds.
 
 \subsection{Cross ratio of fixed points.}
   Supposing that $f$ and $g$ are randomly chosen hyperbolic elements we want to discuss the probability of their axes crossing.  If $f$ has fixed points $z_1,z_2$ and $g$ has fixed points $w_1,w_2$. We identified the formula for the fixed points above at (\ref{fp}) and if we notate the random variables (matrix entries) $a,c$ for $f$ and $\alpha,\beta$ for $g$ we have
 \begin{eqnarray*} 
z_1,z_2 & = & \frac{1}{\bar c}\left( i \Im m(a) \pm \sqrt{\Re e(a)^2-1}\right), \hskip15pt |a|^2=1+|c|^2\\
w_1,w_2 & = & \frac{1}{\bar \beta}\left( i \Im m(\alpha) \pm \sqrt{\Re e(\alpha)^2-1}\right), \hskip15pt |\alpha|^2=1+|\beta|^2
\end{eqnarray*}
and as both elements are hyperbolic we have $\Re e(a)\geq1$ and $\Re e(\alpha)\geq 1$.
We put $U= i \Im m(a) + \sqrt{\Re e(a)^2-1}$ and $V= i \Im m(\alpha) + \sqrt{\Re e(\alpha)^2-1}$
Then
\begin{eqnarray*} [z_1,w_1,z_2,w_2]  
& = & \frac{4 \sqrt{\Re e(a)^2-1} \sqrt{\Re e(\alpha)^2-1}}{\,\bar c \, \bar\beta\,\left(\frac{U}{\bar c} -\frac{V}{\bar \beta} \right) \left(\frac{-\bar U}{\bar c}-\frac{-\bar V}{\bar \beta}\right)} \\
& = &  \frac{4 \sqrt{\Re e(a)^2-1} \sqrt{\Re e(\alpha)^2-1}}{2 \Re e [U\bar V] -c\bar \beta \frac{|U|^2}{|c|^2}-\bar c \beta \frac{|V|^2}{|\beta|^2} } =  \frac{2 \sqrt{\Re e(a)^2-1} \sqrt{\Re e(\alpha)^2-1}}{ \Re e [U\bar V] - \Re e[c\bar \beta] }  
\end{eqnarray*} 
as we recall $1=|z_i|=|U|/|c|$ and similarly $|V|/|\beta|=1$.  Thus we want to understand the statistics of the cross ratio,  and in particular to determine when
\begin{equation}
[z_1,w_1,z_2,w_2]  =  \frac{2 \sqrt{\Re e(a)^2-1} \sqrt{\Re e(\alpha)^2-1}}{ \Re e [U\bar V] - \Re e[c\bar \beta] }   \geq 1 
\end{equation}
We have
\begin{eqnarray*} a=\frac{1}{\sin\theta}e^{i\phi}, \;\; \theta\in_u[0,\pi/2], \phi\in_u[0,2\pi], &&
c=\cot \theta e^{i\delta}, \;\; \delta \in_u[0,2\pi] \\
\alpha =\frac{1}{\sin \eta }e^{i\psi}, \;\; \eta \in_u[0,\pi/2], \psi\in_u[0,2\pi] &&
\beta =\cot \eta e^{i\zeta}, \;\; \zeta \in_u[0,2\pi] 
\end{eqnarray*}
Then $ \sqrt{\Re e(a)^2-1} = \sqrt{\frac{\cos^2\phi}{\sin^2\theta}-1}, \hskip10pt \sqrt{\Re e(\alpha)^2-1} =\sqrt{\frac{\cos^2\psi}{\sin^2\eta}-1}$, 
 $\Phi=\arg c\bar\beta$ is uniformly distributed in $[0,2\pi]$  and 
 \[  \Re e [U\bar V] - \Re e[c\bar \beta] = \frac{\sin \phi}{\sin\theta}\; \frac{\sin \psi}{\sin\eta} +\sqrt{\frac{\cos^2\phi}{\sin^2\theta}-1}\sqrt{\frac{\cos^2\psi}{\sin^2\eta}-1} -\cot \eta\cot \theta \cos \Phi \]
 This gives
\begin{eqnarray*}
\lefteqn{ \frac{2 \sqrt{\Re e(a)^2-1} \sqrt{\Re e(\alpha)^2-1}}{ \Re e [U\bar V] - \Re e[c\bar \beta] }  }\\ & = &\frac{2 \sqrt{ \cos^2\phi -\sin^2\theta}\sqrt{\cos^2\psi-\sin^2\eta }}{ \sin \phi \;  \sin \psi +\sqrt{ \cos^2\phi -\sin^2\theta}\sqrt{\cos^2\psi-\sin^2\eta }-\cos \eta\cos \theta \cos \Phi} \\
& = &\frac{2 \sqrt{ 1-X^2 }\sqrt{1-Y^2}}{ XY+\sqrt{ 1-X^2 }\sqrt{1-Y^2} -  \cos \Phi}  = Z
\end{eqnarray*}
where we define the random variables
\[ X =\frac{\sin \phi}{\cos \theta}, \;\;\;\; {\rm and} \;\;\;\; Y= \frac{\sin \psi}{ \cos\eta} \]
In order for $Z\geq 1$ we need $|X|\leq 1$, $|Y|\leq 1$ and 
\begin{equation} \sqrt{ 1-X^2 }\sqrt{1-Y^2}    \geq \cos \Phi - XY\end{equation}
  If this last condition holds,  then $ [z_1,w_1,z_2,w_2] \geq 1$ requires 
\begin{equation}   \sqrt{ 1-X^2 }\sqrt{1-Y^2}  \geq XY - \cos \Phi 
\end{equation} 
Notice that $X$,  $Y$ and $\Phi\in_u[0,2\pi]$ are independent,  with $X$ and $Y$ identically distributed.  Unfortunately $\sqrt{ 1-X^2 }\sqrt{1-Y^2} \pm  XY $ is difficult to find directly as $\sqrt{ 1-X^2 }\sqrt{1-Y^2}$ and $XY$ are not independent.  We therefore write
\[ X = \sin S, \;\;\;\ S\in [-\frac{\pi}{2},\frac{\pi}{2}],  \hskip10pt    Y = \sin T, \;\;\;\ T\in [-\frac{\pi}{2},\frac{\pi}{2}] \]
so that 
 \[ \sqrt{ 1-X^2 }\sqrt{1-Y^2} \pm XY  = \cos(S\mp T) \]
and we have the two requirements
\begin{equation}\label{2conds} \cos(S\mp T) \geq \pm \cos(\Phi) \end{equation}
Following the arguments of \S 5 we have the p.d.fs 
\begin{eqnarray*}
X & {\rm with \;\; p.d.f.} &F_X(x) = \frac{2}{\pi^2 x} \log \Big| \frac{1+x}{1-x} \Big|, \hskip10pt -1\leq x \leq 1 \\
S& {\rm with \;\; p.d.f.} & F_S(\theta) = \frac{2}{\pi^2} \cot(\theta) \log \Big| \frac{1+\sin(\theta)}{1-\sin(\theta)} \Big|, \hskip10pt -\frac{\pi}{2}\leq \theta \leq \frac{\pi}{2} 
\end{eqnarray*}
We can remove various symmetries and redundancies for the situation to simplify.  For instance we may assume $S\geq 0$ and reduce to ranges where $\cos$ is either increasing or decreasing so we can remove it.  We quickly come to the following conditions equivalent to (\ref{2conds}) with $S$ and $T$ identically distributed as above and $\Phi\in_u[0,\pi/2]$,
\[ 0\leq S, \;\;\;\;-\Phi \leq S-T \leq \Phi, \;\;\;\; {\rm and}\;\;\;\; S+T+\Phi \leq \pi \]
This now sets up an integral which we implemented on Mathematica numerically and which returned the value $0.429\ldots$.  We also ran an experiment using random numbers generated by Mathematica to construct the associated matrices
\[ A=\left(\begin{array}{cc} a & c \\ \bar c & \bar a  \end{array} \right),  \;\;\;B=\left(\begin{array}{cc} \alpha & \beta\\
\bar \beta & \bar \alpha \end{array} \right),  \]
where $a=\frac{e^{i\theta_1}}{\sin(\eta_1)}$, $c=\cot(\eta_1)\; e^{i\theta_2}$, $\alpha=\frac{e^{i\psi_1}}{\sin(\eta_2)}$, $\beta=\cot(\eta_2)\; e^{i\psi_2}$ and distributed 
 \[ \theta_1,\theta_2,\psi_1,\psi_2 \in_u[0,2\pi], \;\;\;\; \eta_1,\eta_2\in_u[0,\pi/2] \]
 We put $\gamma = \gamma(A,B)=\tr [A,B]-2$.

\scalebox{0.275}{\includegraphics*[viewport=250 -100 1800 800]{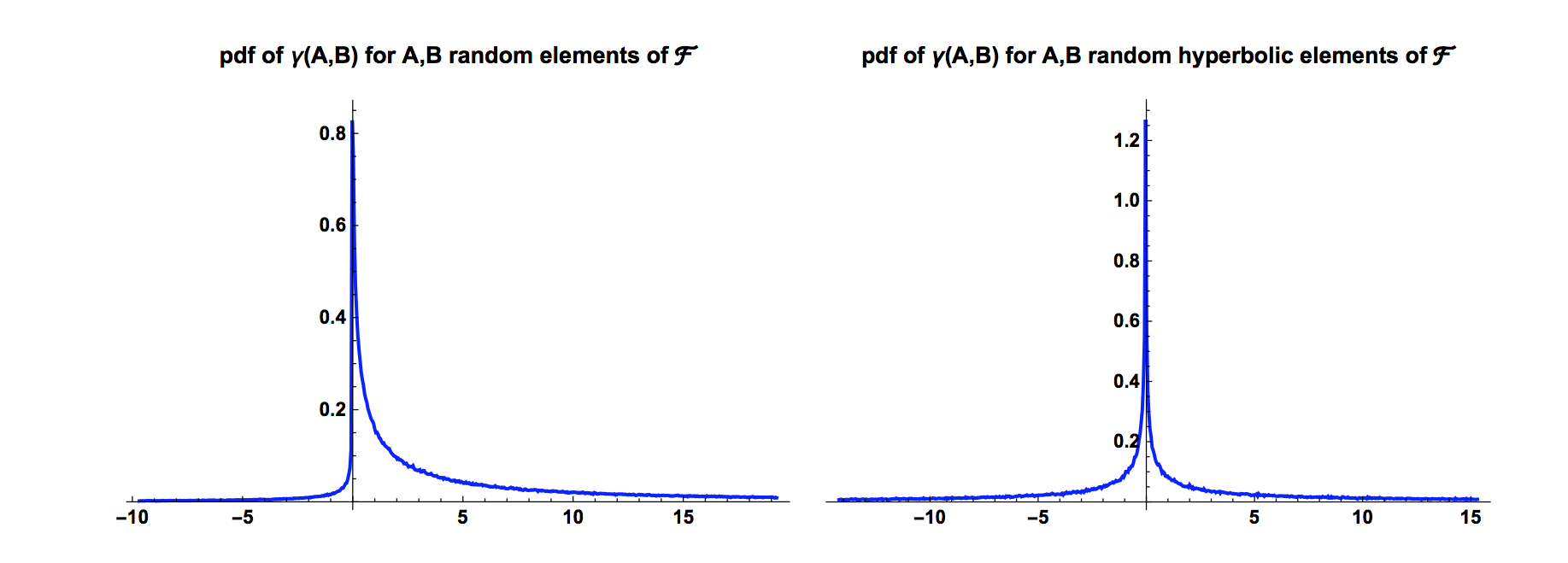}} \\
Left: Histogram of $\gamma(A,B)$ values. \\ Right: Histogram of $\gamma(A,B)$ values conditioned by $A$ and $B$ hyperbolic. 

\medskip

We ran through about $10^7$ random matrix pairs of hyperbolic generators and found the probability that $\gamma<0$ to be about $0.429601$.  

\begin{theorem}  Let $f,g$ be randomly chosen hyperbolic elements of ${\cal F}$.  Then the probability that the axes of $f$ and $g$ cross is $\approx 0.429$. 
\end{theorem}

In contrast,  we have the following theorem.  

\begin{theorem}  Let $\zeta_1,\zeta_1$ and $\eta_1,\eta_2$ be two pairs of points,  each randomly and uniformly chosen on the circle.  Let $\alpha$ be the hyperbolic line between $\zeta_1$ and $\zeta_2$ and $\beta$ the hyperbolic line between $\eta_1$ and $\eta_2$. Then the probability that $\alpha$ and $\beta$ cross is $\frac{1}{3}$.
\end{theorem}
\noindent{\bf Proof.}   We can forget the points come in pairs and label them $z_i$, $i=1,2,3,4$ in order around the circle. There are three different cases all with the same probability.
{\bf 1.} $z_1$ connects to $z_2$,  hence $z_3$ to $z_4$ and the lines are disjoint.
{\bf 2.}  $z_1$ connects to $z_3$,  hence $z_2$ to $z_4$ and the lines intersect.
{\bf 3.} $z_1$ connects to $z_4$,  hence $z_2$ to $z_3$ and the lines are disjoint.

The result now follows.  \hfill $\Box$

\medskip

Together these theorems quantify the degree to which the fixed points are correlated on the circle.  However what we would like to understand is the probability 
\[ \Pr\{\gamma<0| f,g \; \mbox{hyperbolic and $\langle f,g\rangle$ is discrete} \}. \]
Notice that $\gamma(A,B)\in [-4,0]$ implies $\tr^2[A,B]-4\in[-4,0]$ and $[A,B]$ is elliptic and of finite order on a countable subset of $[-4,0]$.
\begin{corollary} If $f,g\in{\cal F}$ are randomly chosen and if $\gamma(f,g)\in [-4,0]$,  then $\langle f,g\rangle$ is almost surely not discrete.
\end{corollary}

We found that $0.266818$ ($\frac{4}{15}$ ?) of our $10^7$ pairs of hyperbolic elements had $\gamma<-4$ while $0.162394$ of the pairs had $-4<\gamma<0$ and so were not discrete and free with probability one.  About $\frac{1}{9}$ of our pairs failed J\o rgensen's test for discreteness, \cite{Jorg}.

\medskip

\scalebox{0.333}{\includegraphics*[viewport=100 0 1000 750]{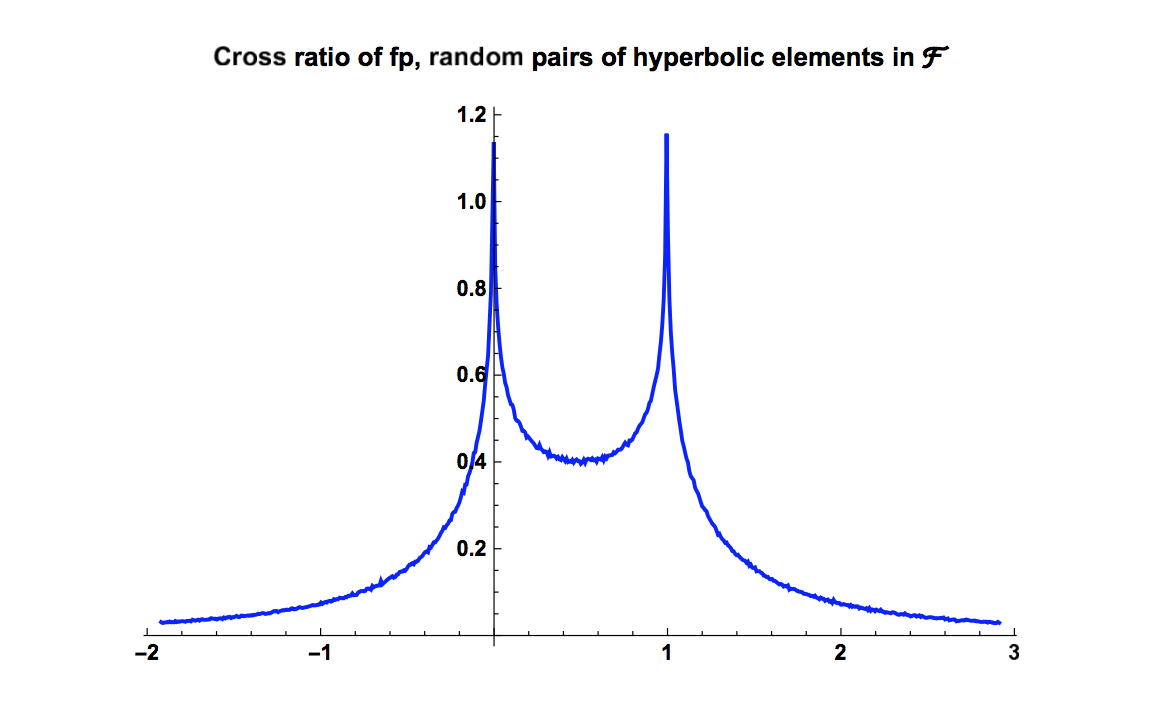}} \\
Histogram of the cross ratio of the fixed points of a randomly chosen pair of hyperbolic elements.

\medskip

In the histogram above the singularities are at $0$ and $1$.  We make the observation that it seems quite likely that ${\rm Pr}\{[z_1,w_1,z_2,w_2]\geq 1\}=\frac{1}{5}$.  It is somewhat of a chore to calculate the cross ratio distribution $X_{cr}$ of four randomly selected point on the circle.  This is done in \cite{GJM} and the distribution is similar to that above,  with singularities at $0$ and $1$.  However for that distribution the probability that 
\[ {\rm Pr}\{X_{cr}<0\}={\rm Pr}\{0<X_{cr}<1\}={\rm Pr}\{X_{cr>1}\}=\frac{1}{3} \]
(as can be seen from the action of the group $S_4$ on the cross ratio,  \cite{Beardon}). This shows the distributions are definitely different.

We next turn to a discussion of positive results for discreteness.

 \section{Discreteness}  
 
 We now have an easy lower bound for the probability a group generated by two random elements of ${\cal F}$ is discrete based on the following Klein combination theorem (or ``ping pong'' lemma).
 
 \begin{lemma}\label{pingpong}  Let $f_i$ $i=1,2,\ldots, n$ be hyperbolic transformations of the disk whose isometric disks are all disjoint.  Then the group generated by these hyperbolic transformations  $\langle f_1,f_2,\ldots,f_n \rangle$ is discrete and isomorphic to the free group $F_n$.
 \end{lemma}
 
We have already seen that the probability that the isometric disks of a randomly chosen  $f\in {\cal F}$ are disjoint is $\frac{1}{2}$.  We can slightly generalise this using Corollary \ref{cor2.6}.

\begin{lemma}  Let $\alpha$ and $\beta$ be arcs on $\IS^1$ with uniformly randomly chosen midpoints $\zeta_\alpha$ and $\zeta_\beta$ and subtending angles $\theta_\alpha$ and $\theta_\beta$ uniformly chosen from $[0,\pi]$.  The the probablity that $\alpha$ and $\beta$ meet is $\frac{1}{2}$.
\end{lemma}

\noindent{\bf Proof.}  The smaller arc subtended between $\zeta_\alpha$ and $\zeta_\beta$ has length $\Theta = \arg(\zeta_\alpha \overline{\zeta_\beta})$ and is uniformly distributed in $[0,\pi]$.  Then $\alpha$ and $\beta$ are disjoint if $\Theta-\theta_\alpha/2-\theta_\beta/2 \geq 0$.  Since Corollary \ref{cor2.6} tells us that $2\Theta - \theta_\alpha-\theta_\beta$ is uniformly distributed in $[-2\pi,2\pi]$ the probability this number is positive is $\frac{1}{2}$. \hfill $\Box$.

\medskip

Using Lemma \ref{pingpong} this quickly gives us the obvious bound that if $f,g\in {\cal F}$ are randomly chosen, then the probability that $\langle f,g \rangle$ is discrete is at least $\frac{1}{64}$.  For $n$ generator groups this number is at least $2^{-(2n-1)!}$.  However we are going to have to build a bit more theory to prove the following substantial improvements of these estimates.

 \begin{theorem} \label{thma} The probability that randomly chosen $f,g\in {\cal F}$ generate a discrete group $\langle f,g\rangle$ is at least $\frac{1}{20}$.
 \end{theorem}
 
 \begin{theorem} \label{thmb}   The probability that two randomly chosen hyperbolic transformation $f,g\in {\cal F}$ generate a discrete group $\langle f,g\rangle$ is at least $\frac{1}{5}$.
 \end{theorem}
 
  \begin{theorem}  Let $f,g$ be randomly chosen parabolic elements in ${\cal F}$.  Then the probability $\langle f,g\rangle$ is discrete is at least $\frac{1}{6}$,
  \[ \Pr\{ \langle f,g \rangle \;\; \mbox{is discrete given $ f,g\in {\cal F}$ are parabolic}  \} \geq  \frac{1}{6} \]
  \end{theorem}
  
  Notice that $f$ is parabolic or the identity if and only if $\Re e(a)\in \{\pm1\}$. Theorem \ref{thma} follows from Theorem \ref{thmb} and the fact that the probablity we choose two hyperbolic elements is independent and of probablity equal to $\frac{1}{4}$.
 
\section{Random arcs on a circle.}

Let $\alpha$ be an arc on the circle $\IS$.  We denote its midpoint by $m_\alpha\in \IS$ and its arclength by $\ell_\alpha\in [0,2\pi]$.  Conversely,  given $m_\alpha\in \IS$ and $\ell_\alpha\in [0,2\pi]$ we determine a unique arc $\alpha = \alpha(m_\alpha,\ell_\alpha)$ with this data.  

A random arc $\alpha$ is the arc uniquely determined when we choose $m_\alpha\in \IS$ uniformly (equivalently $\arg(m_\alpha)\in_u[0,2\pi])$ and $\ell_\alpha \in_u[0,2\pi]$.  We will abuse notation and also refer to random arcs when we restrict to $\ell_\alpha\in_u [0,\pi]$ as for the case of isometric disk intersections.  We will make the distinction clear in context.

A simple consequence of our earlier result is the following corollary.

\begin{corollary}  If $m_\alpha,m_\beta \in_u\IS$ and $\ell_\alpha,\ell_\beta\in_u[0,\pi]$,  then
\[ \Pr \{\alpha \cap \beta = \emptyset\} = \frac{1}{2} \]
\end{corollary}
We need to observe the following lemma.

\begin{lemma} \label{parabolic} If $m_\alpha,m_\beta \in_u\IS$ and $\ell_\alpha,\ell_\beta\in_u[0,2\pi]$,  then
\[ \Pr \{\alpha \cap \beta = \emptyset\} = \frac{1}{6} \]
\end{lemma}
\noindent{\bf Proof.}  We need to calculate the probability that the argument of $\zeta = m_\alpha \overline{m_\beta}$ is greater than $(\ell_\alpha+\ell_\beta)/2$.  Now $\theta = \arg(\zeta)$ is uniformly distributed in $[0,\pi]$.  The joint distribution is uniform,  and so we calculate
\begin{eqnarray*} \Pr\{\theta \geq \ell_\alpha+\ell_\beta \} & = &  \frac{1}{\pi^3} \int\int\int_{\{\theta\geq \alpha + \beta\}} 1 \; d\theta\,d\alpha\, d\beta \\  
& = &  \frac{1}{\pi^3}  \int_{0}^{\pi} \int_{0}^{\theta} \int_{0}^{\theta-\alpha}    d\, \beta   \,d\alpha \,  d\theta = \frac{1}{6} 
\end{eqnarray*}
and the result follows. \hfill $\Box$
  
  \medskip
  Next we consider the probablity of disjoint pairs of arcs.

\begin{lemma}  Let $m_{\alpha_1},m_{\alpha_2},m_{\beta_1},m_{\beta_2} \in_u\IS$ and $\ell_\alpha,\ell_\beta\in_u[0,\pi]$.  Set
\[ \alpha_i=\alpha(m_{\alpha_i},\ell_{\alpha_i}), \hskip10pt \beta_i=\alpha(m_{\beta_i},\ell_{\beta_i}) \]
Then the probability that all the arcs $\alpha_i,\beta_i$, $i=1,2$ are disjoint is $1/20$, 
\[ \Pr \{(\alpha_1\cap\alpha_2)\cup(\beta_1\cap\beta_2)\cup (\alpha_1\cap\beta_1)  \cup (\alpha_1\cap\beta_2)  \cup (\alpha_2\cap\beta_1)  \cup (\alpha_2\cap\beta_2)  = \emptyset\} =  \frac{1}{20} \]
\end{lemma}
 \noindent{\bf Proof.}
 We first observe that the events 
 \[(\alpha_1\cap\beta_1)  =\emptyset, \;\;  (\alpha_1\cap\beta_2)=\emptyset, \;\;     (\alpha_2\cap\beta_1) =\emptyset, \;\;    (\alpha_2\cap\beta_2) =\emptyset \]
 are not independent as (for among other reasons) $\alpha_1$ and $\alpha_2$,  and similarly $\beta_1$ and $\beta_2$ may overlap.  The probability that $(\alpha_2\cap\beta_2) =\emptyset$ and $(\alpha_2\cap\beta_2) =\emptyset$ we have already determined to be equal to $\frac{1}{4}=\frac{1}{2}\times\frac{1}{2}$.  The result now follows from the next lemma. \hfill $\Box$
 
 \begin{lemma}\label{1/5}  Let $m_{\alpha_1},m_{\alpha_2},m_{\beta_1},m_{\beta_2} \in_u\IS$ and $\ell_\alpha,\ell_\beta\in_u[0,\pi]$.  Set
\[ \alpha_i=\alpha(m_{\alpha_i},\ell_{\alpha_i}), \hskip10pt \beta_i=\alpha(m_{\beta_i},\ell_{\beta_i}) \]
and suppose we are given that $(\alpha_1\cap\alpha_2)=(\beta_1\cap\beta_2)=\emptyset$.  Then the probability that all the arcs $\alpha_i$ are disjoint from the arcs $\beta_j$, $i,j=1,2$  is $1/5$, 
\[ \Pr \{(\alpha_1\cap\beta_1)  \cup (\alpha_1\cap\beta_2)  \cup (\alpha_2\cap\beta_1)  \cup (\alpha_2\cap\beta_2)  = \emptyset\} =  \frac{1}{5} \]
\end{lemma}
\noindent{\bf Proof.}  Conditioned by the assumption that $\alpha_1$ and $\alpha_2$ are disjoint,  and that $\beta_1$ and $\beta_2$ are disjoint,  we have note the events 
\[(\alpha_1\cap\beta_1)  =\emptyset, \;\;  (\alpha_1\cap\beta_2)=\emptyset, \;\;     (\alpha_2\cap\beta_1) =\emptyset, \;\;    (\alpha_2\cap\beta_2) =\emptyset \]
are independent.  A little trigonometry reveals that
\[ \alpha_i \cap \beta_j =\emptyset \leftrightarrow  \frac{\ell_{\alpha}+\ell_{\beta}}{2} \leq   2\arcsin \frac{|m_{\alpha_i}-m_{\beta_j}|}{2} = \arg (m_{\alpha_i}\overline{m_{\beta_j}}) \]
Now the four variables $\theta_{i,j} = \arg (m_{\alpha_i}\overline{m_{\beta_j}})$, $i,j=1,2$,  are uniformly distributed in $[0,\pi]$ and independent.  We are requiring
\[ \min_{i,j}  \theta_{i,j} \geq   \frac{\ell_{\alpha}+\ell_{\beta}}{2} \]
Now $\frac{\ell_{\alpha}+\ell_{\beta}}{2}=\psi$ is uniformly distributed in $[0,\pi]$ and 
\begin{equation} \Pr\{ \min_{i,j}  \theta_{i,j} \geq \psi \} = (1-\frac{\psi}{\pi})^4 \label{df} \end{equation}
Since  
\begin{equation} \frac{1}{\pi} \; \int_{0}^{\pi} (1-\frac{\psi}{\pi})^4 =\frac{1}{5} \end{equation} 
the result claimed follows.   \hfill $\Box$

\medskip
In passing we further note that equation (\ref{df})  gives us a density function $\rho(\psi) = 4 (1-\frac{\psi}{\pi})^3$ and hence an expected value of 
\begin{eqnarray*}
\frac{4}{\pi^2} \int_{0}^{\pi} \psi (1-\frac{\psi}{\pi})^3 \; d\psi  & = & 4 \int_{0}^{1} (1-t) t^3\; dt =  \frac{1}{5}.
\end{eqnarray*}

 \medskip
 
 Generalising this result for a greater number of disjoint pairs  of arcs quickly gets quite complicated.  We state without proof given here the following which we will not use.
\begin{lemma}  Let $m_{\alpha_1},m_{\alpha_2},m_{\beta_1},m_{\beta_2},m_{\gamma_1},m_{\gamma_2} \in_u\IS$ and $\ell_\alpha,\ell_\beta,\ell_\gamma\in_u[0,\pi]$.  Set
\[ \alpha_i=\alpha(m_{\alpha_i},\ell_{\alpha_i}), \hskip10pt \beta_i=\alpha(m_{\beta_i},\ell_{\beta_i}), \hskip10pt \gamma_i=\alpha(m_{\gamma_i},\ell_{\gamma_i}) \]
Then the probability that all the arcs $\alpha_i,\beta_i,\gamma_i$, $i=1,2$ are all disjoint is $\frac{3}{1000}$.  
\end{lemma}

One can get results if there is additional symmetry.  For instance if the lengths of all the arcs are the same.

\begin{theorem}  Let $m_{i_1},m_{i_2} \in_u\IS^1$, $i=2,\ldots,n$ and $\ell_\alpha\in_u[0,\pi]$.  Then the probability that the arcs $\alpha_{ij}=\alpha(m_{i_j},\ell_\alpha)$ are disjoint is
\begin{equation} \frac{1}{(2n) n!} \int_{0}^{1} \sum_{k=0}^{[2-x]} (-1)^k \left(\begin{array}{c} n \\ k\end{array} \right) (2-x-k)^{n} \; dx \end{equation}
\end{theorem}
\noindent{\bf Proof.} We cyclically order the set $\{m_{i_i}:i=2,\ldots, n, j=1,2\}$ and let $\theta_k$ be the angle between the $k^{th}$ and $k+1^{st}$ point (mod $k$).  Then $\sum_{k=1}^{2n} \theta_k = 2\pi$.  The arcs are disjoint if $\theta_k \geq \ell_\alpha$.  First we have $2n-1$ independent random variables $\{\theta_k\}_{k=1}^{2n-1}$ whose minimum must exceed $\alpha$,  and second they also must satisfy $2\pi - \sum_{k=1}^{2n-1} \theta_k \geq \ell_\alpha$.  The first gives us a factor $\frac{1}{2n}$,  and for the second we note that the sum of $m$ uniformly distributed random variables in $[0,1]$ has the Irwin-Hall distribution, 
\begin{equation}\label{IH}
F_n(x) = \frac{1}{(m-1)!} \sum_{k=0}^{[x]} (-1)^k \left(\begin{array}{c} m \\ k\end{array} \right) (x-k)^{m-1}
\end{equation}
Thus 
\[ \Pr\left\{2-\frac{\ell_\alpha}{\pi} \geq \sum_{k=1}^{2n-1} \frac{\theta_k}{\pi}\right\} = \int_{0}^{2-t} F_{2n-1} (t)\; dt  \]
The result follows. \hfill $\Box$

\medskip
As an example,  for two pairs of equi-length arcs we have
\[ F_3(x) = \left\{ \begin{array}{ll}  x^2/2, & 0\leq x\leq 1 \\  (-2x^2+6x-3)/2, & 1\leq x \leq 2  \\  (x^2-6x+9)/2, & 2\leq x \leq 3 \end{array}\right. \]
We see that
\begin{eqnarray*}
\int_{0}^{2-t} F_3(x)\; dx & = &  \int_{0}^{1} F_3(x)\; dx + \int_{1}^{2-t} F_3(x) \; dx = \frac{1}{6} +  \frac{2}{3}-\frac{t}{2}-\frac{t^2}{2}+\frac{t^3}{3}  \\
\int_{0}^{1} \int_{0}^{2-t} F_3(x) \,dx dt & = &   \frac{1}{6} + \int_{0}^{1} \frac{2}{3}-\frac{t}{2}-\frac{t^2}{2}+\frac{t^3}{3} \, dt = \frac{1}{6}+\frac{1}{3} = \frac{1}{2} 
\end{eqnarray*}
and so the probability that two pairs of random equi-arclength arcs with arclength uniformly distributed in $[0,\pi]$,  are disjoint is $\frac{1}{8}$.  Similarly for three pairs the probability is $\frac{9}{200}$.
  
 \section{Random arcs to M\"obius groups.}
 
 Given data $m_{\alpha_1},m_{\alpha_2} \in \IS$ with arclength $\ell_\alpha\in [0,\pi]$ we see,  just as  above, that the arcs centered on the $m_{\alpha_i}$ and on length $\ell_\alpha$ determine a matrix which can be calculated by examination of the isometric circles.  We have 
\begin{equation} \label{Adef} A= \left(\begin{array}{cc} a & c\\ \bar c & \bar a \end{array} \right), \;\;\;  c = i \sqrt{m_{\alpha_1}\, m_{\alpha_2}}\;  \cot \frac{\ell_\alpha}{2},\;\;\;\; a=
  i \sqrt{ \overline{m_{\alpha_1}} \,m_{\alpha_2}}\;  {\rm cosec} \frac{\ell_\alpha}{2} 
  \end{equation}
 where we make a consistent choice of sign by ensuring 
 \[  \frac{c}{a} = m_{\alpha_1} \cos \frac{\ell_\alpha}{2}  \]
 Of course interchanging $m_{\alpha_1}$ and $m_{\alpha_2}$ sends $a$ to $-\bar a$,  and so the data actually uniquely determines the cyclic group $\langle f \rangle$ generated by the 
associated M\"obius transformation 
 \[ f(z) = - m_{\alpha_2 }\; \frac{z+m_{\alpha_1}\cos \frac{\ell_\alpha}{2}}{z\, \cos \frac{\ell_\alpha}{2}+m_{\alpha_1}} \]
 and not necessarily $f$ itself.  
 
 As a consequence we have the following theorem.
 \begin{theorem}  There is a one-to-one correspondence between collections of $n$ pairs of random arcs and $n$-generator Fuchsian groups.  A randomly chosen $\langle f \rangle\subset {\cal F}$ corresponds uniquely to $m_{\alpha_1}, m_{\alpha_2} \in_u \IS^1$ and $\ell_\alpha\in_u [0,\pi]$.
 \end{theorem}
 Notice also that if we recognise the association of cyclic groups with the data and say two cyclic groups are close if they have close generators,  then this association is continuous.  
 
 \medskip
 
 We have already seen that for a pair of hyperbolic elements if all the isometric disks are disjoint then the ``ping ping'' lemma implies discreteness of the groups in question.  Then the association between Fuchsian groups and random arcs  quickly establishes Theorems \ref{thma} and \ref{thmb} via Lemma \ref{1/5}.
 
 \medskip

If $f$ is a parabolic element of ${\cal F}$,  then the isometric circles are adjacent and meet at the fixed point.  Conversely,  if two random arcs of arclength $\ell_\alpha$ are adjacent we have $\arg(m_{\alpha_1}\overline{m_{\alpha_2}})=\ell_\alpha$,  and from (\ref{Adef})
 \[ a=
  i (\cos\frac{\ell_\alpha}{2}  +i\sin \frac{\ell_\alpha}{2}  ) {\rm cosec} \frac{\ell_\alpha}{2}  = -1+ i \cot\frac{\ell_\alpha}{2} \]
 and $\tr^2(A)-4=0$ so that $A$ represents a parabolic transformation.  Similarly if the arcs overlap,  then $\tr^2(A)\leq 2$ and $A$ represents an elliptic transformation.
 
 \begin{theorem}  Let $f,g$ be randomly chosen parabolic elements in ${\cal F}$.  Then the probability $\langle f,g\rangle$ is discrete is at least $\frac{1}{6}$. 
  \end{theorem}
  \noindent{\bf Proof.}  As $f$ and $g$ are parabolic,  their isometric disks are tangent and the point of intersection lies in a random arc of arclength uniformly distributed in $[0,2\pi]$.  Discreteness follows from the ``ping pong''  lemma and  Lemma \ref{parabolic}.  \hfill $\Box$

\end{document}